\def\titlep{Classification and realizations of type III factor representations
of Cuntz-Krieger algebras associated with quasi-free states}
\documentclass[11pt]{article}
\usepackage{graphicx,ifthen}
\usepackage{amssymb}
\usepackage{amsmath}

 at12pt

\newcommand{\qed}{\hbox{\rule[-2pt]{3pt}{6pt}}}
\newcommand{\qedh}{\hfill\qed \\}

\newcommand{\vep}{\varepsilon}

\newtheorem{Thm}{Theorem}[section]

\newtheorem{ex}[Thm]{Example}

\newtheorem{lem}[Thm]{Lemma}

\newtheorem{prop}[Thm]{Proposition}

\newtheorem{cor}[Thm]{Corollary}

\def\cal#1{\mathcal #1}
\def\con{{\cal O}_{n}}

\def\pr{{\it Proof.}\quad}

\def\coa{{\cal O}_{A}}
\def\co#1{{\cal O}_{#1}}
\def\ltn{l_{2}({\bf N})}

\def\brl{branching law}
\def\bfsnl{{\rm BFS}_{N}(\Lambda)}

\addtocounter{footnote}{1}
\def\sftt#1{
\setcounter{equation}{0}
\addtocounter{footnote}{1}
\section{#1}
}

\def\ssft#1{\subsection{#1}}

\begin{document}
%
%
\def\autherp{Katsunori Kawamura}
\def\emailp{e-mail: kawamura@kurims.kyoto-u.ac.jp.}
\def\addressp{{\small {\it College of Science and Engineering Ritsumeikan University,}}\\
{\small {\it 1-1-1 Noji Higashi, Kusatsu, Shiga 525-8577, Japan}}
}

\def\infw{\Lambda^{\frac{\infty}{2}}V}
\def\zhalfs{{\bf Z}+\frac{1}{2}}
\def\ems{\emptyset}
\def\pmvac{|{\rm vac}\!\!>\!\! _{\pm}}
\def\vac{|{\rm vac}\rangle _{+}}
\def\dvac{|{\rm vac}\rangle _{-}}
\def\ovac{|0\rangle}
\def\tovac{|\tilde{0}\rangle}
\def\expt#1{\langle #1\rangle}
\def\zph{{\bf Z}_{+/2}}
\def\zmh{{\bf Z}_{-/2}}
\def\brl{branching law}
\def\bfsnl{{\rm BFS}_{N}(\Lambda)}
\def\scm#1{S({\bf C}^{N})^{\otimes #1}}
\def\mqb{\{(M_{i},q_{i},B_{i})\}_{i=1}^{N}}
\def\zhalf{\mbox{${\bf Z}+\frac{1}{2}$}}
\def\zmha{\mbox{${\bf Z}_{\leq 0}-\frac{1}{2}$}}
\newcommand{\mline}{\noindent
\thicklines
\setlength{\unitlength}{.1mm}
\begin{picture}(1000,5)
\put(0,0){\line(1,0){1250}}
\end{picture}
\par
 }
\def\ptimes{\otimes_{\varphi}}
\def\delp{\Delta_{\varphi}}
\def\delps{\Delta_{\varphi^{*}}}
\def\gamp{\Gamma_{\varphi}}
\def\gamps{\Gamma_{\varphi^{*}}}
\def\sem{{\sf M}}
\def\hdelp{\hat{\Delta}_{\varphi}}
\def\tilco#1{\tilde{\co{#1}}}
\def\ndm#1{{\bf M}_{#1}(\{0,1\})}
\def\cdm#1{{\cal M}_{#1}(\{0,1\})}
\def\tndm#1{\tilde{{\bf M}}_{#1}(\{0,1\})}
\def\sck{{\sf CK}_{*}}
\def\hdel{\hat{\Delta}}
\def\ba{\mbox{\boldmath$a$}}
\def\bb{\mbox{\boldmath$b$}}
\def\bc{\mbox{\boldmath$c$}}
\def\be{\mbox{\boldmath$e$}}
\def\bp{\mbox{\boldmath$p$}}
\def\bq{\mbox{\boldmath$q$}}
\def\bu{\mbox{\boldmath$u$}}
\def\bv{\mbox{\boldmath$v$}}
\def\bw{\mbox{\boldmath$w$}}
\def\bx{\mbox{\boldmath$x$}}
\def\by{\mbox{\boldmath$y$}}
\def\bz{\mbox{\boldmath$z$}}
\def\bomega{\mbox{\boldmath$\omega$}}
\def\N{{\bf N}}
\def\lxm{L_{2}(X,\mu)}
\def\rnl{RN_{loc}(X)}
\def\scn{S({\bf C}^{n})}
\def\bonk{\{1,\ldots,n\}^{k}}
\def\bfset#1#2{{\rm BFS}_{#1}(#2)}
%
%
%
\setcounter{section}{0}
\setcounter{footnote}{0}
\setcounter{page}{1}
\pagestyle{plain}

%
%
\title{\titlep}
\author{\autherp\thanks{\emailp}
\\
\addressp}
\date{}
\maketitle

%
%
\begin{abstract}
We completely classify type III factor representations
of Cuntz-Krieger algebras associated with quasi-free states
up to unitary equivalence.
Furthermore,
we realize these representations on concrete Hilbert spaces
without using GNS construction.
Free groups and their type ${\rm II}_{1}$ factor representations
are used in these realizations.
\end{abstract}


%
%
\sftt{Introduction}
\label{section:first}
Representations of operator algebras
are ingredients necessary for quantum field theory \cite{Araki,Haag}.
We have studied representations of C$^{*}$-algebras.
Any $*$-representation of any separable C$^{*}$-algebra
on a separable Hilbert space
is canonically decomposed into direct integrals of factor representations,
and a factor representation is either type I or II or III (\cite{Blackadar2006}, $\S$ III.5).
For any simple Cuntz-Krieger algebra $\coa$ \cite{CK},
there exists neither type ${\rm I}_{n}$ ($1\leq n<\infty)$ 
nor type ${\rm II}_{1}$ nondegenerate representation
because $\coa$ is purely infinite (\cite{RS}, R\o rdam, Proposition 4.4.2) and 
a purely infinite C$^{*}$-algebra has no nondegenerate 
lower semicontinuous trace (\cite{Blackadar2006}, Proposition V.2.2.29).
Type I, especially, irreducible representations were studied in \cite{CKR02}.
In \cite{Okayasu}, 
type III's were constructed as Gel'fand-Na\v{\i}mark-Segal (=GNS) representations from
quasi-free states \cite{EFW,EL} and
every isomorphism classes of von Neumann algebras $N_{\pi}$ generated by
such representations $\pi$ were completely classified.
However, 
even if $N_{\pi}$ and $N_{\pi^{'}}$ are isomorphic,
$\pi$ and $\pi^{'}$ are not equivalent in general.

In this paper,
we completely classify type III factor representations
of Cuntz-Krieger algebras associated with quasi-free states
up to unitary equivalence.
Next, we realize these representations on
concrete Hilbert spaces without using GNS construction from states.
Free groups and their type ${\rm II}_{1}$ factor representations
are used in these realizations.

%
%
\ssft{Type III factor representations of Cuntz-Krieger algebras
associated with quasi-free states}
\label{subsection:firstone}
In this subsection,
we review states and representations of Cuntz-Krieger algebras
by \cite{EL,Okayasu} according to notations and symbols in \cite{TS07}.
We state that a matrix $A\in M_{n}({\bf C})$ is 
{\it nondegenerate} if any column and any row are not zero;
$A$ is {\it irreducible} 
if for any $i,j\in\{1,\ldots,n\}$,
there exists $k\in {\bf N}\equiv \{1,2,3,\ldots\}$ such that
$(A^{k})_{i,j}\ne 0$ where $A^{k}=A\cdots A$ ($k$ times).
For $2\leq n<\infty$,
let $\cdm{n}$ denote the set of all 
irreducible nondegenerate $n\times n$ matrices with entries $0$ or $1$,
which is not a permutation matrix.
Define
%
%
\begin{equation}
\label{eqn:nondegenerate}
\cdm{*}\equiv \cup \{\cdm{n}:n\geq 2\}.
\end{equation}
Let 
$I_{0}^{n}\equiv \{(x_{i})_{i=1}^{n}\in {\bf R}^{n}:\mbox{for all }i,\,0<x_{i}<1\}$.
For $\ba=(a_{i})_{i=1}^{n}\in I_{0}^{n}$,
define $\hat{\ba}\equiv {\rm diag}(a_{1},\ldots,a_{n})\in M_{n}({\bf R})$.
For $A\in \cdm{n}$,
let $\hat{\ba}A$ denote the product of matrices $\hat{\ba}$ and $A$.
Define the set $\Lambda(A)$ of vectors by
%
%
\begin{equation}
\label{eqn:pfe}
\Lambda(A)\equiv \{\ba\in I_{0}^{n}:PFE(\hat{\ba}A)=1\}
\end{equation}
where $PFE(X)$ denotes the Perron-Frobenius eigenvalue 
of an irreducible non-negative matrix $X$ \cite{BP,Seneta}. 
For $A\in \cdm{n}$ and $\ba=(a_{i})_{i=1}^{n}\in \Lambda(A)$,
let $(x_{i})_{i=1}^{n}$ 
denote the Perron-Frobenius eigenvector 
of $\hat{\ba}A$ such that $x_{1}+\cdots+x_{n}=1$ and
let $s_{1},\ldots,s_{n}$ denote the canonical generators of $\coa$.
Define the state $\rho_{\ba}$ over $\coa$ by
%
%
\begin{equation}
\label{eqn:kmsthree}
\rho_{\ba}(s_{J}s_{K}^{*})
=\delta_{JK}a_{j_{1}}\ldots a_{j_{m-1}}x_{j_{m}}
\end{equation}
when $s_{J}s_{K}^{*}\ne 0$
for $J=(j_{1},\ldots,j_{m})\in\{1,\ldots,n\}^{m}$
and $K\in \cup_{l\geq 1}\{1,\ldots,n\}^{l}$
where $s_{J}=s_{j_{1}}\cdots s_{j_{m}}$.
The state $\rho_{\ba}$ is called {\it quasi-free} \cite{Okayasu}.
The relation between the original style of $\rho_{\ba}$ in \cite{Okayasu}
and (\ref{eqn:kmsthree}) is shown in Lemma 3.1 of \cite{TS07}.
Then  the following holds.
%
%
\begin{Thm}
\label{Thm:typetwo}
For $\ba\in\Lambda(A)$,
let $\varpi_{\ba}$ denote the GNS representation of $\coa$ by $\rho_{\ba}$.
Then the von Neumann algebra $M_{\ba}\equiv \varpi_{\ba}(\coa)^{''}$ is an 
approximately finite dimensional factor of type ${\rm III}$.
Furthermore,
Connes' classification of the type of $M_{\ba}$ \cite{Connes}
is given as follows:
\begin{enumerate}
\item
If there exist $p_{1},\ldots,p_{n}\in {\bf N}$ and $0<\lambda<1$
such that  the greatest common divisor of the set $\{p_{1},\ldots,p_{n}\}$
is $1$ and $\ba=(\lambda^{p_{1}},\ldots,\lambda^{p_{n}})$,
then $p_{1},\ldots,p_{n},\lambda$ are uniquely determined by $\ba$,
and $M_{\ba}$ is of type ${\rm III}_{\lambda}$.
\item
If the assumptions in (i) do not hold,
then 
$M_{\ba}$ is of type ${\rm III}_{1}$.
\end{enumerate}
\end{Thm}
Theorem \ref{Thm:typetwo} is a reformulation of Theorem 4.2 in \cite{Okayasu}
without use of terminology of KMS states.

From Theorem \ref{Thm:typetwo},
$0<\lambda\leq 1$ is uniquely determined from a given $\ba\in \Lambda(A)$
such that $\varpi_{\ba}$ is of type ${\rm III}_{\lambda}$.
We write this as $\lambda(\ba)$.
Then the family $\{\varpi_{\ba}:\ba\in\Lambda(A)\}$ of 
type III factor representations of $\coa$
are roughly classified by real numbers $\lambda(\ba)$. 

%
%
\ssft{Classification}
\label{subsection:firsttwo}
Let $\cdm{*}$ be as in (\ref{eqn:nondegenerate}) and 
let $A\in\cdm{*}$.
For $\Lambda(A)$ in (\ref{eqn:pfe}), 
$\ba\in\Lambda(A)$ and $\varpi_{\ba}$ in Theorem \ref{Thm:typetwo},
we give the complete classification of 
unitary equivalence classes of $\{\varpi_{\ba}:\ba\in\Lambda(A)\}$
as follows.
%
%
\begin{Thm}
\label{Thm:mainone}
For $\ba,\bb\in \Lambda(A)$,
$\varpi_{\ba}$ and $\varpi_{\bb}$ are unitarily equivalent 
if and only if $\ba=\bb$.
\end{Thm}

\noindent
Immediately, we see that 
for two quasi-free states $\rho$ and $\rho^{'}$ over $\coa$,
their GNS representations $\pi_{\rho}$ and $\pi_{\rho^{'}}$
are unitarily equivalent if and only if $\rho=\rho^{'}$.
Remark that $\varpi_{\ba}$ and $\varpi_{\bb}$ are unitarily equivalent 
if and only if they are quasi-equivalent 
because they are of type III and their representation
spaces are separable (\cite{Dixmier}, $\S$ 5.6.6).
Therefore $\varpi_{\ba}$ and $\varpi_{\bb}$ are quasi-equivalent 
if and only if $\ba=\bb$.
From Theorem \ref{Thm:mainone}, the following holds.
%
%
\begin{cor}
\label{cor:equivalence}
The set $\{[\varpi_{\ba}]:\ba\in\Lambda(A)\}$ 
of unitary equivalence classes of representations
and $\Lambda(A)$ are in one-to-one correspondence. 
\end{cor}

%
%
\ssft{Realizations}
\label{subsection:firstthree}
We realize $\varpi_{\ba}$ in Theorem \ref{Thm:typetwo} in this subsection.
Let $L_{2}[0,1]$ denote the Hilbert space of all 
complex-valued square integrable functions on the closed interval $[0,1]$.
At the beginning, we construct a representation of $\coa$ on $L_{2}[0,1]$.
For $A=(A_{ij})\in\cdm{n}$ and $\ba=(a_{i})_{i=1}^{n}\in\Lambda(A)$,
let $(x_{i})_{i=1}^{n}$ denote the Perron-Frobenius eigenvector of 
$\hat{\ba}A$  such that $x_{1}+\cdots+x_{n}=1$.
Define real numbers $\{c_{i}\}_{i=0}^{n}$ and $\{b_{ij}\}_{i,j=1}^{n}$ by
\[c_{0}\equiv 0,\quad c_{i}\equiv \sum_{k=1}^{i}x_{k},\quad
b_{ij}\equiv c_{i-1}-a_{i}\sum_{k=1}^{j}(1-A_{ik})x_{k}\quad(i,j=1,\ldots,n).\]
Then $0=c_{0}<c_{1}<\cdots<c_{n}=1$ 
and $b_{i1}\geq \cdots \geq b_{in}$ for each $i$.
Define closed intervals $R_{i}\equiv [c_{i-1},\,c_{i}]$ and
$V_{ij}\equiv [a_{i}c_{j-1}+b_{ij},\,a_{i}c_{j}+b_{ij}]$ for $i,j=1,\ldots,n$.
Then $\bigcup_{j=1;A_{ij}=1}^{n}V_{ij}=R_{i}$
and $V_{ij}\cap V_{ij^{'}}$ is a null set when $j\ne j^{'}$.
Define the representation $\eta_{\ba}$ of $\coa$ on $L_{2}[0,1]$ by
%
%
\begin{equation}
\label{eqn:repone}
\{\eta_{\ba}(s_{i})\phi\}(y)\equiv 
\left\{
\begin{array}{ll}
a_{i}^{-1/2}\phi((y-b_{ij})/a_{i})\quad &(y\in V_{ij},\,A_{ij}=1),\\
\\
0\quad &(\mbox{otherwise})
\end{array}
\right.
\end{equation}
for $i=1,\ldots,n$, $\phi\in L_{2}[0,1]$ and $y\in [0,1]$.

Next,
let ${\bf F}_{n}$ denote the free group
with generators $\xi_{1},\ldots,\xi_{n}$ and
let $U$ denote the left regular representation of ${\bf F}_{n}$
on $l_{2}({\bf F}_{n})$.

By using $(L_{2}[0,1],\eta_{\ba})$ and $(l_{2}({\bf F}_{n}),U)$,
we define a new representation of $\coa$ as follows.
%
%
\begin{Thm}
\label{Thm:exist}
Define the representation
$\Pi_{\ba}$ of $\coa$ on the Hilbert space
$L_{2}[0,1]\otimes l_{2}({\bf F}_{n})$ by
%
%
\begin{equation}
\label{eqn:largepi}
\Pi_{\ba}(s_{i})\equiv \eta_{\ba}(s_{i})\otimes U_{\xi_{i}}
\quad(i=1,\ldots,n).
\end{equation}
Let $\{e_{g}:g\in {\bf F}_{n}\}$ denote the standard basis 
of $l_{2}({\bf F}_{n})$ and
define the cyclic subspace ${\cal K}_{\ba}$ of 
$L_{2}[0,1]\otimes l_{2}({\bf F}_{n})$ by
%
%
\begin{equation}
\label{eqn:largek}
{\cal K}_{\ba}\equiv \overline{\Pi_{\ba}(\coa)({\bf 1}\otimes e_{\vep})}
\end{equation}
where ${\bf 1}$ denotes the constant function on $[0,1]$ with value $1$ and 
$\vep$ denotes the unit of ${\bf F}_{n}$.
Then $({\cal K}_{\ba},\Pi_{\ba}|_{{\cal K}_{\ba}})$ 
is a type {\rm III} factor representation of $\coa$
which is unitarily equivalent to $\varpi_{\ba}$
in Theorem \ref{Thm:typetwo}.
\end{Thm}

From Theorem \ref{Thm:mainone} and \ref{Thm:exist}, the following is verified
for $\eta_{\ba}$ in (\ref{eqn:repone}).
%
%
\begin{prop}
\label{prop:eta}
For each $\ba,\bb\in\Lambda(A)$,
$\eta_{\ba}$ and $\eta_{\bb}$ are unitarily equivalent 
if and only if $\ba=\bb$.
\end{prop}
From Theorem \ref{Thm:mainone} and Proposition \ref{prop:eta},
the classification of $\{\varpi_{\ba}:\ba\in\Lambda(A)\}$
is equivalent to that of $\{\eta_{\ba}:\ba\in\Lambda(A)\}$.
The idea of the construction of $\eta_{\ba}$ has originated as 
representations of Cuntz-Krieger algebras 
arising from interval dynamical systems \cite{KS1,CKR01,PFO01}.
When $A_{ij}=1$ for all $i,j$, $\coa\cong \con$ and 
it is known that
$\eta_{\ba}$ is irreducible for each $\ba\in\Lambda(A)$
(Appendix \ref{section:appone}).

In $\S$ \ref{section:second}, we will prove Theorem \ref{Thm:mainone} 
and \ref{Thm:exist}.
In $\S$ \ref{section:third}, we will show examples.

%
%
\sftt{Proofs of main theorems}
\label{section:second}
We prove main theorems in this section.
%
%
\ssft{Proof of Theorem \ref{Thm:mainone}}
\label{subsection:secondone}
In order to prove Theorem \ref{Thm:mainone},
we prepare two lemmata.
%
%
\begin{lem}
\label{lem:inequality}
For $A=(A_{ij})\in\cdm{n}$,
let $\ba=(a_{i})_{i=1}^{n},\bb=(b_{i})_{i=1}^{n}\in\Lambda(A)$ and 
let $\bx=(x_{i})_{i=1}^{n}$ and $\by=(y_{i})_{i=1}^{n}$
be Perron-Frobenius eigenvectors of $\hat{\ba}A$ 
and $\hat{\bb}A$, respectively such that
$x_{1}+\cdots+x_{n}=1$ and $y_{1}+\cdots+y_{n}=1$.
For $i,j=1,\ldots,n$,
define
%
%
\begin{equation}
\label{eqn:dij}
D_{ij}\equiv \sqrt{a_{i}b_{i}A_{ij}}.
\end{equation}
Assume $\ba\ne \bb$.
\begin{enumerate}
\item
There exists $0<c<1$ such that the following holds for any $i=1,\ldots,n$:
%
%
\begin{equation}
\label{eqn:sh}
\sum_{j=1}^{n}D_{ij}\sqrt{x_{j}y_{j}}\leq c\sqrt{x_{i}y_{i}}.
\end{equation}
\item
For $m\geq 2$,
define the positive real number $T_{m}$ by
%
%
\begin{equation}
\label{eqn:tm}
T_{m}\equiv \sum_{(j_{1},\ldots,j_{m})
\in\{1,\ldots,n\}^{m}}D_{j_{1}j_{2}}\cdots D_{j_{m-1}j_{m}}
\sqrt{x_{j_{m}}y_{j_{m}}}.
\end{equation}
Then $T_{m}\to 0$ when $m\to\infty$.
\end{enumerate}
\end{lem}
%
%
\pr
(i)
Fix $i\in \{1,\ldots,n\}$. Let $\bv\equiv(\sqrt{a_{i}A_{ij}x_{j}})_{j=1}^{n}$
and $\bw\equiv (\sqrt{b_{i}A_{ij}y_{j}})_{j=1}^{n}$.
From the Schwarz inequality
and the assumption of $\ba,\bb,\bx,\by$,
%
%
\begin{equation}
\label{eqn:vw}
\sum_{j=1}^{n}D_{ij}\sqrt{x_{j}y_{j}}
=\langle \bv|\bw\rangle\leq \|\bv\|\|\bw\|=\sqrt{x_{i}}\sqrt{y_{i}}.
\end{equation}
Define $c_{i}\equiv \sum_{j=1}^{n}D_{ij}\sqrt{x_{j}y_{j}}/\sqrt{x_{i}y_{i}}$.
Then $0\leq c_{i}\leq 1$.
Since the L.H.S. of (\ref{eqn:sh}) is not zero,  $c_{i}>0$.
Assume $c_{i}=1$.
Then we see that there exists $k>0$ such that $\bv=k\bw$.
This implies $\bx=\by$. From this, $\ba=\bb$.
This is a contradiction. Hence $c_{i}<1$.
Let $c\equiv \max_{i}\{c_{i}\}$.
Then $0<c<1$ and we obtain that 
$\sum_{j=1}^{n}D_{ij}\sqrt{x_{j}y_{j}}=c_{i}\sqrt{x_{i}y_{i}}
\leq c\sqrt{x_{i}y_{i}}$ for each $i$.
From this, the statement holds.

\noindent
(ii)
Let $c$ be as in (i).
From (i), $T_{m}\leq cT_{m-1}$ when $m\geq 3$.
Let $\sqrt{\bx}\equiv (\sqrt{x_{i}})_{i=1}^{n}$
and $\sqrt{\by}\equiv (\sqrt{y_{i}})_{i=1}^{n}$.
Then $\|\sqrt{\bx}\|=\|\sqrt{\by}\|=1$.
From (i) and the Schwarz inequality,
%
%
\begin{equation}
\label{eqn:ttwo}
T_{2}=\sum_{i,j=1}^{n}D_{ij}\sqrt{x_{j}y_{j}}\leq \sum_{i=1}^{n}c\sqrt{x_{i}y_{i}}
\leq c\langle\sqrt{\bx}|\sqrt{\by}\rangle \leq c.
\end{equation}
From these, we see that $T_{m}\leq c^{m-1}$ for each $m\geq 2$.
This implies the statement.
\qedh
%
%
\begin{lem}
\label{lem:orth}
Let $A,\ba,\bb,\bx,\by$ be as in Lemma \ref{lem:inequality}.
Assume that $\coa$ acts on a Hilbert space ${\cal H}$.
Let $\xi,\zeta$ be unit vectors in ${\cal H}$ 
and let $\rho$ and $\rho^{'}$ denote vector states over $\coa$ 
with respect to $\xi$ and $\zeta$, respectively.
When $\ba\ne \bb$, the following holds.
\begin{enumerate}
\item 
Assume that 
$\rho$ and $\rho^{'}$ satisfy
%
%
\begin{equation}
\label{eqn:tthree}
\rho(s_{J}s_{J}^{*})=a_{j_{1}}\cdots a_{j_{m-1}}x_{j_{m}},\quad 
\rho^{'}(s_{J}s_{J}^{*})=b_{j_{1}}\cdots b_{j_{m-1}}y_{j_{m}}
\end{equation}
when $s_{J}\ne 0$
for each $J=(j_{1},\ldots,j_{m})\in\{1,\ldots,n\}^{m}$, $m\geq 1$.
Then $\langle\xi|\zeta\rangle=0$.
\item 
Assume that 
there exist $i_{0},i_{0}^{'}\in\{1,\ldots,n\}$ such that
$\rho$ and $\rho^{'}$ satisfy
\[\rho(s_{J}s_{J}^{*})=\kappa_{i_{0},j_{1}}
a_{j_{1}}\cdots a_{j_{m-1}}x_{j_{m}},\quad
\rho^{'}(s_{J}s_{J}^{*})=
\kappa^{'}_{i_{0}^{'},j_{1}}b_{j_{1}}\cdots b_{j_{m-1}}y_{j_{m}}\]
when $s_{J}\ne 0$
for each $J=(j_{1},\ldots,j_{m})\in\{1,\ldots,n\}^{m}$, $m\geq 1$
where $\kappa_{i_{0},j_{1}}\equiv a_{i_{0}}A_{i_{0},j_{1}}/x_{i_{0}}$
and $\kappa^{'}_{i_{0}^{'},j_{1}}
\equiv b_{i_{0}^{'}}A_{i_{0}^{'},j_{1}}/y_{i_{0}^{'}}$.
Then $\langle\xi|\zeta\rangle=0$.
\end{enumerate}
\end{lem}
%
%
\pr
(i)
Let $T_{m}$ be as in (\ref{eqn:tm}) and $E_{J}=s_{J}s_{J}^{*}$.
Since $\sum_{J\in\{1,\ldots,n\}^{m}}E_{J}=I$ and 
$|\langle\xi|E_{J}\zeta\rangle|\leq \sqrt{\rho(E_{J})}\sqrt{\rho^{'}(E_{J})}$,
%
%
\begin{equation}
\label{eqn:norm}
|\langle\xi|\zeta\rangle|\leq 
\sum_{J\in\{1,\ldots,n\}^{m}}\sqrt{\rho(E_{J})}\sqrt{\rho^{'}(E_{J})}=T_{m}.
\end{equation}
This holds for each $m\geq 2$.
From Lemma \ref{lem:inequality}(ii), the statement holds.

\noindent
(ii)
Let $D_{ij}$ be as in (\ref{eqn:dij}).
From the proof of (i),
\[
\begin{array}{rl}
|\langle\xi|\zeta\rangle|\leq &
\sum_{J\in\{1,\ldots,n\}^{m}}\sqrt{\rho(E_{J})}\sqrt{\rho^{'}(E_{J})}\\
=&
q\sum_{(j_{1},\ldots,j_{m})\in\{1,\ldots,n\}^{m}}
A_{i_{0},j_{1}}A_{i_{0}^{'},j_{1}}D_{j_{1}j_{2}}\cdots D_{j_{m-1}j_{m}}\sqrt{x_{j_{m}}y_{j_{m}}}\\
\leq & qT_{m}\to 0\qquad(m\to\infty)
\end{array}
\]
where $q\equiv\{a_{i_{0}}b_{i_{0}^{'}}/(x_{i_{0}}y_{i_{0}^{'}})\}^{1/2}$.
Hence the statement holds.
\qedh

\noindent
{\it Proof of Theorem \ref{Thm:mainone}.}
It is sufficient to show 
that 
$\varpi_{\ba}$ and $\varpi_{\bb}$ are not unitarily equivalent 
if $\ba\ne \bb$.
Assume that 
$\varpi_{\ba}$ and $\varpi_{\bb}$ are unitarily equivalent when $\ba\ne \bb$.
We can assume that they act on the same Hilbert space ${\cal H}$
such that both $\varpi_{\ba}(s_{i})$ and $\varpi_{\bb}(s_{i})$ are
written as the same operator $s_{i}$ on ${\cal H}$ for each $i$.
Then there exist two cyclic unit vectors $\Omega_{\ba}$ and $\Omega_{\bb}$
such that the vector states over $\coa$ associated 
with $\Omega_{\ba}$ and $\Omega_{\bb}$
coincide with $\rho_{\ba}$ and $\rho_{\bb}$, respectively.
From Lemma \ref{lem:orth}(i), $\langle\Omega_{\ba}|\Omega_{\bb}\rangle=0$.

For $J,K\in\{1,\ldots,n\}^{*}$ such that $s_{J}s_{K}^{*}\ne 0$,
let $X_{JK}\equiv \langle\Omega_{\ba}|s_{J}s_{K}^{*}\Omega_{\bb}\rangle$.
Define
$\xi\equiv d_{1}s_{J}^{*}\Omega_{\ba}$ 
and $\zeta\equiv d_{2}s_{K}^{*}\Omega_{\bb}$
where $d_{1}$ and $d_{2}$ are normalization constants
of $\xi$ and $\zeta$, respectively.
Then we see that $\xi$ and $\zeta$ satisfy the assumption in 
Lemma \ref{lem:orth}(ii).
Hence $X_{JK}=(d_{1}d_{2})^{-1}\langle\xi|\zeta\rangle=0$
for any $J,K$.
Since ${\rm Lin}\langle \{s_{J}s_{K}^{*}\Omega_{\bb}:J,K\}\rangle$
is dense in $\overline{\coa\Omega_{\bb}}$,
$\Omega_{\ba}=0$.
This is a contradiction.
Therefore $\varpi_{\ba}$ and $\varpi_{\bb}$
are not unitarily equivalent.
\qedh

%
%
\ssft{Proof of Theorem \ref{Thm:exist}}
\label{subsection:secondtwo}

\noindent
{\it Proof of Theorem \ref{Thm:exist}.}
Define the state $\rho$ over $\coa$ by
$\rho\equiv \langle{\bf 1}|\eta_{\ba}(\cdot){\bf 1}\rangle$.
By the definition of $\eta_{\ba}$, we can verify that
%
%
\begin{equation}
\label{eqn:kms}
\rho(s_{J}s_{J}^{*})=a_{j_{1}}\cdots a_{j_{m-1}}x_{j_{m}}
\end{equation}
when $s_{J}\ne 0$
for each $J=(j_{1},\ldots,j_{m})\in \{1,\ldots,n\}^{m}$.
Define the state $\tilde{\rho}$ over $\coa$ by
%
%
\begin{equation}
\label{eqn:tilde}
\tilde{\rho}
\equiv \langle {\bf 1}
\otimes e_{\vep}|\Pi_{\ba}(\cdot)({\bf 1}\otimes e_{\vep})\rangle.
\end{equation}
By using (\ref{eqn:kms}),
we can verify that $\tilde{\rho}=\rho_{\ba}$.
From Theorem \ref{Thm:typetwo}
and the uniqueness of GNS representation
up to unitary equivalence, the statement holds.
\qedh

Especially, when $A_{ij}=1$ for all $i,j$,
$\coa\cong \con$ and
the state $\rho$ in (\ref{eqn:kms}) satisfies
%
%
\begin{equation}
\label{eqn:rhotwo}
\rho(s_{J}s_{K}^{*})=\sqrt{a_{J}a_{K}}
\quad(J,K\in\bigcup_{l\geq 1}\{1,\ldots,n\}^{l})
\end{equation}
where $a_{J}\equiv a_{j_{1}}\cdots a_{j_{m}}$ when $J=(j_{1},\ldots,j_{m})$.
On the other hand, the state $\tilde{\rho}$ in (\ref{eqn:tilde}) satisfies
%
%
\begin{equation}
\label{eqn:rhothree}
\tilde{\rho}(s_{J}s_{K}^{*})=\delta_{JK}
a_{J}\quad(J,K\in\bigcup_{l\geq 1}\{1,\ldots,n\}^{l}).
\end{equation}
In order to construct $\tilde{\rho}$ from $\rho$,
we add the Kronecker delta of (\ref{eqn:rhothree}) to $\rho$
by using the left regular representation of the free group ${\bf F}_{n}$.

%
%
\sftt{Examples}
\label{section:third}
We show examples in this section.
%
%
\ssft{$\Lambda(A)$}
\label{subsection:thirdone}
For $n\geq 2$,
let $A=(A_{ij})\in\cdm{n}$.
We show examples of $\Lambda(A)$ in (\ref{eqn:pfe}).
%
%
\begin{ex}
\label{ex:one}
{\rm
If $A_{ij}=1$ for each $i,j$,
then the following holds (\cite{Izumi}, $\S$ 2.1):
\[\Lambda(A)=\{(a_{i})_{i=1}^{n}:a_{1}+\cdots+a_{n}=1,\,a_{i}>0\mbox{ for all }i\}.\]
}
\end{ex}
%
%
\begin{ex}
\label{ex:two}
{\rm
If $A=\left(
\begin{array}{cc}
1&1\\
1&0\\
\end{array}
\right)$,
then
\[\Lambda(A)=\left\{\left(x,\,\frac{1}{x}-1\right): \frac{1}{2}<x<1\right\}.\]
}
\end{ex}
%
%
\begin{ex}
\label{ex:three}
{\rm
If $A=\left(
\begin{array}{ccc}
0&1&1\\
1&0&1\\
1&1&1\\
\end{array}
\right)$,
then
\[\Lambda(A)=\left\{\left(x,\,y,\,\frac{1-xy}{1+x+y+xy}\right): 
0<x,y<1\right\}.\]
For example,
$\ba=(1/3,1/3,1/2)$ belongs to $\Lambda(A)$.
In this case, $(1/4,1/4,1/2)$
is the Perron-Frobenius eigenvector of $\hat{\ba}A$.
The representation $(L_{2}[0,1],\eta_{\ba})$ of $\coa$ in (\ref{eqn:repone}) 
is given as follows:
\[
\left\{
\begin{array}{rl}
\{\eta_{\ba}(s_{1})\phi\}(y)=&\sqrt{3}\chi_{R_{1}}(y)\phi(3y+1/4),\\
\\
\{\eta_{\ba}(s_{2})\phi\}(y)=&
\sqrt{3}\chi_{W_{1}}(y)\phi(3y-3/4)+\sqrt{3}\chi_{W_{2}}(y)\phi(3y-1/2),\\
\\
\{\eta_{\ba}(s_{3})\phi\}(y)=&\sqrt{2}\chi_{R_{3}}(y)\phi(2y-1)
\end{array}
\right.
\]
for $\phi\in L_{2}[0,1]$ and $y\in [0,1]$
where $R_{1}\equiv [0,1/4]$, $R_{3}\equiv[1/2,1]$,
$W_{1}\equiv [1/4,1/3]$ and $W_{2}\equiv [1/3,1/2]$, and
$\chi_{Y}$ denotes the characteristic function of $Y\subset [0,1]$.
}
\end{ex}

%
%
\ssft{Several realizations of representations of $\con$}
\label{subsection:thirdtwo}
%
%
\begin{ex}
\label{ex:four}
{\rm
Assume $A_{ij}=1$ for each $i,j$.
For $\ba=(a_{i})_{i=1}^{n}\in\Lambda(A)$,
let $c_{0}\equiv 0$, $c_{i}\equiv \sum_{j=1}^{i}a_{j}$ for $i=1,\ldots,n$.
Then  $\eta_{\ba}$ in (\ref{eqn:repone})
is given as follows:
%
%
\begin{equation}
\label{eqn:etatwo}
\{\eta_{\ba}(s_{i})\phi\}(y)
\equiv a_{i}^{-1/2}\chi_{[c_{i-1},c_{i}]}(y)\phi((y-c_{i-1})/a_{i})
\quad(i=1,\ldots,n)
\end{equation}
for $y\in [0,1]$ and $\phi\in L_{2}[0,1]$.
}
\end{ex}
%
%
\begin{ex}
\label{ex:five}
{\rm
Let $\{\xi_{1},\xi_{2}\}$ and $\vep$ denote generators 
and the unit of ${\bf F}_{2}$, respectively.
Assume that $A=(A_{ij})\in \cdm{2}$ and $A_{ij}=1$ for each $i,j=1,2$.
From Example \ref{ex:one},
we see that $\Lambda(A)=\{(a,b):a+b=1,\,a,b>0\}$.
For $\ba=(a,b)\in\Lambda(A)$,
we introduce a representation $\Pi_{\ba}^{'}$ of $\co{2}$  
on the Hilbert space $l_{2}({\bf N}\times {\bf F}_{2})$
which is unitarily equivalent to $\Pi_{\ba}$ in 
Theorem \ref{Thm:exist}.
Define the representation $\eta_{\ba}^{'}$ of $\co{2}$ on $\ltn$ by
%
%
\begin{equation}
\label{eqn:dash}
\eta_{\ba}^{'}(s_{1})e_{n}\equiv \sqrt{a}e_{2n-1}-\sqrt{b}e_{2n},\quad 
\eta_{\ba}^{'}(s_{2})e_{n}\equiv\sqrt{b}e_{2n-1}+\sqrt{a}e_{2n}
\end{equation}
for $n\in {\bf N}$
where $\{e_{n}:n\in {\bf N}\}$ denotes the standard basis of $\ltn$.
Since $\eta_{\ba}^{'}(\sqrt{a}s_{1}+\sqrt{b}s_{2})e_{1}=e_{1}$
and $e_{1}$ is cyclic for $(\ltn,\eta_{\ba}^{'})$,
$\eta_{\ba}^{'}$ is
unitarily equivalent to  $\eta_{\ba}$ in (\ref{eqn:repone})
(Appendix \ref{section:appone}).
Define
$\Pi_{\ba}^{'}(s_{i})\equiv \eta^{'}_{\ba}(s_{i})\otimes U_{\xi_{i}}$
for $i=1,2$.
Then $\Pi_{\ba}^{'}$ is a representation of $\co{2}$ on $l_{2}({\bf N}\times {\bf F}_{2})$
which is unitarily equivalent to $\Pi_{\ba}$.
Let $\{e_{n,g}:(n,g)\in {\bf N}\times {\bf F}_{2}\}$
denote the standard basis of $l_{2}({\bf N}\times {\bf F}_{2})$.
Then the following holds:
%
%
\begin{equation}
\label{eqn:otwo}
\left\{
\begin{array}{rl}
\Pi_{\ba}^{'}(s_{1})e_{n,g}=&
\sqrt{a}e_{2n-1,\xi_{1}g}-\sqrt{b}e_{2n,\xi_{1}g},\\
\\
\Pi_{\ba}^{'}(s_{2})e_{n,g}
=&\sqrt{b}e_{2n-1,\xi_{2}g}+\sqrt{a}e_{2n,\xi_{2}g}
\end{array}
\right.
\end{equation}
for each $(n,g)\in {\bf N}\times {\bf F}_{2}$.
Especially, we can verify that 
$\rho_{\ba}=\langle e_{1,\vep}|\Pi_{\ba}^{'}(\cdot)e_{1,\vep}\rangle$.

For $p,q\in {\bf N}$, $\lambda>0$,
assume that $\lambda^{p}+\lambda^{q}=1$ and 
the greatest common divisor of $p$ and $q$ is $1$.
If $\ba=(\lambda^{p},\lambda^{q})$,
then the cyclic subrepresentation of $\Pi_{\ba}^{'}$
with the cyclic vector $e_{1,\vep}$
is a type ${\rm III}_{\lambda}$ factor representation
of $\co{2}$ which is unitarily equivalent to $\varpi_{\ba}$
from Theorem \ref{Thm:typetwo}.
}
\end{ex}


\appendix
\section*{Appendix}
%
%
\sftt{A family of irreducible representations of $\con$}
\label{section:appone}
We show a family of representations of Cuntz algebras \cite{PFO01}
related to $\eta_{\ba}$ in (\ref{eqn:repone}).
Let $S({\bf C}^{n})\equiv\{\bz\in {\bf C}^{n}:\|\bz\|=1\}$.
For $\bz=(z_{i})_{i=1}^{n}\in S({\bf C}^{n})$,
let $GP(\bz)$ denote the class of representations $({\cal H}, \pi)$ with
a cyclic unit vector $\Omega\in {\cal H}$ such that
$\pi(z_{1}s_{1}+\cdots+z_{n}s_{n})\Omega=\Omega$.

We show results of $GP(\bz)$ as follows.
For any $\bz\in S({\bf C}^{n})$,
the class $GP(\bz)$ consists of only one unitary equivalence class.
Hence we can identify the class $GP(\bz)$ and the representative of $GP(\bz)$.
For any $\bz\in S({\bf C}^{n})$, $GP(\bz)$ is irreducible.
For $\bz,\bz^{'}\in S({\bf C}^{n})$, $GP(\bz)$ and $GP(\bz^{'})$
are unitarily equivalent  if and only if $\bz=\bz^{'}$.
For $\bz=(z_{i})_{i=1}^{n}\in S({\bf C}^{n})$, 
$GP(\bz)$ is unitarily equivalent to
the GNS representation by the state $\rho$ over $\con$ defined by
$\rho(s_{J}s_{J^{'}}^{*}) \equiv \overline{z_{J}}z_{J^{'}}$
where $J,J^{'}\in \bigcup_{k\geq 1}\{1,\ldots,n\}^{k}$ and
$z_{J}\equiv z_{j_{1}}\cdots z_{j_{k}}$ for $J=(j_{1},\ldots,j_{k})$.

Especially,
if $\bz=(\sqrt{a_{i}})_{i=1}^{n}$,
then $\eta_{\ba}$ in (\ref{eqn:repone}) is $GP(\bz)$
for $\ba=(a_{i})_{i=1}^{n}\in\Lambda(A)$ with respect to the matrix $A_{ij}=1$ for each $i,j$.

\section*{Acknowledgement}
The author would like to express his sincere thanks to referees 
for their numerous suggestions.
%
%
%

\end{document}